\documentclass[11pt,reqno]{amsart}

\usepackage{amssymb}

\newtheorem{theorem}{Theorem}

\newtheorem{remark}[theorem]{Remark}

\begin{document}

\title[On computing joint invariants of vector fields]{On computing joint invariants
of vector fields}

\author[H. Azad]{H. Azad}

\address{Department of Mathematics and Statistics, King Fahd University,
Saudi Arabia}

\email{hassanaz@kfupm.edu.sa}

\author[I. Biswas]{I. Biswas}

\address{School of Mathematics, Tata Institute of Fundamental Research, Homi
Bhabha Road, Bombay 400005, India}

\email{indranil@math.tifr.res.in}

\author[R. Ghanam]{R. Ghanam}

\address{Virginia Commonwealth University in Qatar, Education City
Doha, Qatar}

\email{raghanam@vcu.edu}

\author[M. T. Mustafa]{M. T. Mustafa}

\address{Department of Mathematics, Statistics and Physics, Qatar
University, Doha, 2713, State of Qatar}

\email{tahir.mustafa@qu.edu.qa}



\date{}

\begin{abstract}
A constructive version of the Frobenius integrability theorem -- that
can be programmed effectively -- is given. This is used in computing
invariants of groups of low ranks and recover examples from a recent
paper of Boyko, Patera and Popoyvich \cite{BPP}.
\end{abstract}

\maketitle

\section{Introduction}

The effective computation of local invariants of Lie algebras of
vector fields is one
of the main technical tools in applications of Lie's symmetry method
to several problems in differential equations -- notably their
classification and explicit solutions of natural equations of
mathematical physics, as shown, e.g., in several papers of Ibragimov
\cite{Ib2}, \cite{Ib3}, and Olver \cite{Ol1}.

The main aim of this paper is to give a constructive procedure that
reduces the determination of joint local invariants of any finite
dimensional Lie algebra of vector fields --- indeed any finite
number of vector fields --- to that of a commuting family
of vector fields. It is thus a constructive version of the
Frobenius integrability theorem -- \cite[p. 422]{Ol1}, \cite[p. 472]{Le},
\cite[p. 92--94]{Ib1} -- which can also be programmed effectively. This is actually
valid for any field of scalars. A paper close to this paper is \cite{Pe1}.

We illustrate the main result by computing joint
invariants for groups of low rank as well as examples from Boyko et
al \cite{BPP}, where the authors have used the method of moving frames, \cite{Ol2},
to obtain invariants.

It is stated in \cite{BPP} that solving the first order system of differential
equations is not practicable. However, it is practicable for at least two reasons.
The local joint invariants in any representation of a Lie algebra as an algebra
of vector fields are the same as those of a commuting family of operators.
Moreover, one needs to take only operators that are generators for the full
algebra. For example, if the the Lie algebra is semisimple with Dynkin diagram having
$n$ nodes, then one needs just $2n$ basic operators to determine invariants.

Another reason is that software nowadays can handle symbolic computations very
well.

The main results of the paper are as follows:

\begin{theorem}\label{prop1}
Let $\mathcal L$ be a finite dimensional Lie algebra of
vector fields defined on some open subset $U$ of
$\mathbb{R}^{n}$. Let $X_{1}\, ,\cdots\, ,X_{d}$ be a basis of $\mathcal L$.
Then the following hold:
\begin{enumerate}
\item The algebra of operators whose coefficient matrix is the matrix
of functions obtained from the coefficients of $X_1\, ,\cdots \, ,X_d$ by
reducing it to reduced row echelon form is abelian.

\item The local joint invariants of $\mathcal L$ are the same as those of the
above abelian algebra.
\end{enumerate}
\end{theorem}

\begin{theorem}\label{thm2}
Let $X_1\, , X_2\, ,\cdots\, ,X_d$ be vector fields defined on some open subset of
${\mathbb R}^n$. Then the joint invariants of $X_1\, , X_2\, ,\cdots\, ,X_d$ are
given by the following algorithm:
\begin{enumerate}
\item{[Step 1]} Find the row reduced echelon form of $X_1\, , X_2\, ,\cdots\, ,X_d$, and let
$Y_1\, ,\cdots\, ,Y_r$ be the corresponding vector fields. If this is a commuting
family then stop. Otherwise go to:

\item{[Step 2]} If some $[Y_i\, ,Y_j]\,\not=\, 0$, then set $Y_{r+1}\, :=\,
[Y_i\, ,Y_j]$. Go to Step 1 and substitute $Y_1\, ,\cdots\, ,Y_r\, , Y_{r+1}$
in place of $X_1\, , X_2\, ,\cdots\, ,X_d$.
\end{enumerate}
This process terminates in at most $n$ iterations. If $V_1\, ,\cdots\, ,V_m$ are the
commuting vector fields at the end of the above iterative process, the joint invariants
of $X_1\, , X_2\, ,\cdots\, ,X_d$ coincide with the joint invariants
of $V_1\, ,\cdots\, ,V_m$.
\end{theorem} 

\section{Some examples and proof of Theorem \ref{prop1} and Theorem \ref{thm2}}

Before proving Theorem \ref{prop1}, we give some examples in detail, because
these examples contain all the key ideas of a formal proof and of
computation of local joint invariants of vector fields.

\subsection{Example: The rotations in ${\mathbb R}^3$}

The group ${\rm SO}(3)$ has one
basic invariant in its standard representation, namely $x^2+y^2+z^2$, which
is clear from geometry. Let us recover this by Lie algebra
calculations in a manner that is applicable to all Lie groups.

The fundamental vector fields given by rotations in the coordinate
planes are
$$
I\,=\,y\frac{\partial}{\partial x}-x\frac{\partial}{\partial y}\, , \ J\,=\,
z\frac{\partial}{\partial y}-y\frac{\partial}{\partial z}~ \text{ and } ~
K\,=\, z\frac{\partial}{\partial x}-x\frac{\partial}{\partial z}\, .
$$
The coefficients matrix is
\begin{equation}\label{e1}
\begin{pmatrix}
 y & -x & 0 \\
0 & z & -y \\
z & 0 & -x
\end{pmatrix}\, .
\end{equation}
This is a singular matrix, so its rank is at most two. On the open subset $U$ where
$yz\,\neq \,0$, the rank is two. The rank is two everywhere except at the origin but
we are only interested in the rank on some open set.

The differentiable functions on $U$ simultaneously annihilated by $I\, ,J\, ,K$ are
clearly the same as those of the operators whose coefficient matrix is obtained from
\eqref{e1} reducing it to its row echelon form. Since $I\, ,J$ generate the infinitesimal
rotations, we may delete the last row in \eqref{e1}.
The reduced row echelon form of \eqref{e1} is
$$
\begin{pmatrix}
1 & 0 & \displaystyle \frac{-x}{z} \\
0 & 1 & \displaystyle \frac{-y}{z}
\end{pmatrix}\, .
$$

The operators whose matrix of coefficients is this matrix are
$$ X\,:=\, \frac{\partial}{\partial x}-\frac{x}{z}\frac{\partial}{\partial z}\ \text{ and }
\ Y\,:=\,\frac{\partial}{\partial y}-\frac{y}{z}\frac{\partial}{\partial z}\, .$$
Note that $[X\, ,Y]\,=\,0$. Now, because the fields are commuting, we can compute the
basic invariants of any one of them, say $X$; then $Y$ will operate on the
invariants of $X$.\\

The invariants for $X$ are given by the standard method of Cauchy characteristics
as follows \cite[p. 67]{Ib1}: We want to solve
$$\frac{dx}{1}\,=\,\frac{dy}{0}\,=\,\frac{-zdz}{x}\, .$$
The basic invariants of $X$ are $x^{2}+z^{2}\,=:\,\xi\, ,\ y\,=:\,\eta$.
As $Y$ commutes with $X$, it operates on invariants of $X$.
Now $Y(\xi)\,=\, -2\eta\, ,\ Y(\eta)\,=\,1$. Thus on the invariants of $X$ the field
induced by $Y$ is
$$-2\eta \frac{\partial}{\partial\xi}+\frac{\partial}{\partial \eta}\, .$$
The corresponding characteristic system is $$\frac{d\xi}{-2\eta}\,=\,
\frac{d\eta}{1}\, ,$$ so we get the basic invariant -- which must be a joint
invariant -- as $\xi+\eta^{2}\,=\,x^{2}+y^{2}+z^{2}$.

Examples given below shows what happens if we just work with
finitely many vector fields.

\subsection{Example: The rotations in ${\mathbb R}^n$ with metric signature $(p\, ,q)$,
where $p+q\,=\, n$.} The group $\text{SO}(p,q)$ operates transitively on every
nonzero level set of the
function $\sum_{i=1}^p x^2_i -\sum_{i=1}^q x^2_{i+p}$, and it operates transitively
on the nonzero vectors in the zero level set of this function.
Therefore, it is clear that there is only one
basic joint invariant. Let us recover this by Lie algebra calculations in a manner
that is applicable in general.

The Lie group $\text{SO}(p,q)$ is generated by ordinary rotations in the $(x_1\, ,x_2)$--plane,
the $(x_2\, ,x_3)$--plane, $\cdots \,$, the $(x_{p-1}\, ,x_p)$--plane,
the $(x_{p+1}\, ,x_{p+2})$--plane, $\cdots\, $, the $(x_{p+q-1}\, ,x_{p+q})$--plane,
and hyperbolic rotations in the $(x_p\, ,x_{p+1})$--plane. The fundamental vector
fields generated by these rotations in the coordinate planes are
$$
x_{i+1}\frac{\partial}{\partial{x_i}} - x_{i}
\frac{\partial}{\partial{x_{i+1}}}\, ,i\, \in\,\{1\, ,\cdots\, ,p+q-1\}\setminus
\{p\}\ \text{ and }~\ x_{p+1}\frac{\partial}{\partial{x_p}} +x_{p}
\frac{\partial}{\partial{x_{p+1}}}\, .
$$
The reduced row echelon form is the augmented $(n-1)\times (n-1)$ identity matrix, augmented
by column vector
$$
\frac{x_1}{x_n}\, ,\cdots\, , \frac{x_p}{x_n}\, ,-\frac{x_{p+1}}{x_n}\, ,\cdots\, ,
-\frac{x_{p+q-1}}{x_n}\, .
$$
Thus we get the corresponding vector fields
$$
\frac{\partial}{\partial{x_i}} + \frac{x_{i}}{x_n}\frac{\partial}{\partial{x_n}}
\, ,i\, \leq\, p~\ \text{ and }~\
\frac{\partial}{\partial{x_j}} - \frac{x_j}{x_n}\frac{\partial}{\partial{x_n}}
\,\  ,p\, <\, j\, \leq\, n-1\, .
$$
Since for independent variables $x\, , y\, ,z$,
$$
[\frac{\partial}{\partial{x}} + \frac{x}{z}\frac{\partial}{\partial{z}}\, ,
\frac{\partial}{\partial{y}} + \frac{y\epsilon}{z}\frac{\partial}{\partial{z}}]
\,=\, [\frac{x}{z}\frac{\partial}{\partial{z}}\, , \frac{y\epsilon}{z}
\frac{\partial}{\partial{z}}]\,=\, 0\, ,
$$
where $\epsilon\,=\, \pm 1$, we conclude that these vector fields commute and each such
field operates on the invariants of the remaining. By calculations as in Example 1
we see that the basic joint invariant is $\sum_{i=1}^p x^2_i -\sum_{i=1}^q x^2_{i+p}$.

\subsection{Proof of Theorem \ref{prop1}}

We will use the notation in the statement of Theorem \ref{prop1}. Take a point $p
\,\in\, U$, and let ${\mathcal L}(p)$ be the linear span of $X(p)$ with $X \,\in\,
\mathcal L$. Let $r(p)$ be the
dimension of ${\mathcal L}(p)$, and let $r\,=\,\max \, \{r (p)\}_{p\in U}$. Choose a
point $p$ with $r(p)\,=\,r$.

By renaming the basis for $\mathcal L$, we may assume that $X_{1}(p)\, ,\cdots \,
,X_{r}(p)$ is a basis for ${\mathcal L}(p)$. Therefore, the determinant
$X_1(p)\bigwedge \cdots \bigwedge X_r(p)\,\in\, \bigwedge^r T_p U$ is nonzero.
Hence $X_1(q)\bigwedge \cdots \bigwedge X_r(q)\,\in\, \bigwedge^r T_q U$ is nonzero
for all $q$ in a neighborhood of $p$. In particular, $r(q)\,=\, r(p)\,=\, r$ at all
such points $q$.

Replacing $U$ by this open neighborhood of $p$, we may suppose that $r(q)\,=\,r$ for all
points $q\,\in\, U$. This implies that $X_{r+k}(q)$ is a linear combination of $X_{1}(q)\, ,
\cdots\, ,X_{r}(q)$ with coefficients that depend differentiably on $q\, \in\, U$.
Moreover, for any $X\, ,Y\, \in\, \mathcal L$, as $[X\, ,Y](q)$ is a linear combination
of $X_1(q)\, ,\cdots\, ,X_d(q)$ with scalar coefficients, we see that for
$1\, \leq\, i\, ,j\, \leq\, r$, the Lie bracket $[X_i \, , X_j](q)$ is a linear
combination of $X_1(q)\, ,\cdots\, ,X_r(q)$ with coefficients that depend differentiably
on $q$. Also, for $1\, \leq\, i\, ,j\, \leq\, r$ and any differentiable function $f$,
$$
[X_i\, , fX_j]
$$
is a linear combination of $X_1\, ,\cdots\, ,X_r$ with coefficients that are differentiable
functions. If
$$
X_j\, =\, \sum_{k=1}^n a_{jk} \frac{\partial}{\partial{x_k}}
\, ,~ \ 1\,\leq\, j\, \leq\, r\, ,
$$
we put these operators in reduced row echelon form with coefficients as differentiable
functions. Therefore, taking possibly a smaller open subset of $U$, we obtain a family
of vector fields which span ${\mathcal L}(q)$, $q\, \in\, U$, and is closed under Lie
brackets with differentiable functions as coefficients.
Also, the local invariants for this family are the same as for
$X_{1}\, ,\cdots\, ,X_{r}$.

After changing indices, we may suppose that
$$
X_j\,=\, \frac{\partial}{\partial{x_j}}+\sum_{k=r+1}^n b_{jk}\frac{\partial}{\partial{x_k}}
\, ,~ \ 1\,\leq\, j\, \leq\, r\, .
$$
We want to show that $[X_i\, , X_j]\,=\, 0$ for all $i\, , j\, \leq\, r$,
A straightforward computation shows that
$$
[X_i\, ,X_j]\, \equiv\,0 ~ \ \text{modulo }\ \frac{\partial}{\partial x_{r+1}}\, ,
\cdots \, , \frac{\partial}{\partial{x_n}}\, ,
$$
meaning $[X_i\, ,X_j]\, =\,\sum_{\ell=r+1}^n\phi^\ell_{i,j} \frac{\partial}{\partial
x_\ell}$, where $\phi^\ell_{i,j}$ are smooth functions.
On the other hand, $[X_i\, ,X_j]$ is a linear combination of $X_1\, ,\cdots\, ,X_r$
with functions as coefficients. From these we conclude that $[X_i\, ,X_j]\,
=\, 0$. This completes the proof of the theorem.

\subsection{Proof of Theorem \ref{thm2}}

We use the notation of Theorem \ref{thm2}. Since $Y_1\, ,\cdots\, ,Y_r$ are in row reduced
echelon form, and $[Y_i\, ,Y_j]\,\not=\, 0$, it follows that $[Y_i\, ,Y_j]$ is not
in the linear span of $Y_1\, ,\cdots\, ,Y_r$ with smooth functions as coefficients.
Therefore, when we go back to Step 1 and construct the row
reduced echelon form of $Y_1\, ,\cdots\, ,Y_r\, , [Y_i\, ,Y_j]$, there are $r+1$ vector
fields in it. Consequently, each time we come back and complete Step 1, the number of
vector fields goes up by one. This immediately implies that the process stops after
at most $n$ iterations. The final statement of the theorem is obvious.

\begin{remark} {\rm Theorem \ref{prop1} and Theorem \ref{thm2} are valid 
in algebraic category for any field -- working with the Zariski 
topology. For the field $\mathbb R$, one has the standard refinement 
that $r$ commuting fields of rank $r$ are, in suitable coordinates 
$\frac{\partial}{\partial x_i}$ (Frobenius' theorem). The reason is 
that any nonzero vector field $X$ in suitable local coordinates is
$\frac{\partial}{\partial x_1}$ and any vector field that commutes with
$X$ operates on the invariants of $X$.}
\end{remark}

Let us illustrate Theorem \ref{thm2} by two examples.

Taking the example in \cite[p. 64]{Ol1},
consider the following three vector fields on ${\mathbb R}^3$:
$$
V_+\,:=\, 2y\frac{\partial}{\partial x} + z\frac{\partial}{\partial y}\, ,~
V_0\,:=\, -2x\frac{\partial}{\partial x} + 2z\frac{\partial}{\partial z}\, ,~
V_{-}\,:=\, -x\frac{\partial}{\partial y} + 2y\frac{\partial}{\partial z}\, .
$$
Although they are closed under Lie bracket, we do not need this fact to compute the
joint invariants.

The row reduced echelon form of the matrix of coefficients is
$$
\begin{pmatrix}
1 & 0 & -z/x\\
0 & 1 & 2y/x
\end{pmatrix}
$$
Let
$$
X\, :=\, \frac{\partial}{\partial x} - \frac{z}{x}\frac{\partial}{\partial z}~\
\text{ and }~\ Y\, :=\, \frac{\partial}{\partial y} +
\frac{2y}{x}\frac{\partial}{\partial z}\, .
$$
Sine $[X\, ,Y]\,=\,0$, we stop at this stage. The invariants of $X$ are
$$
\xi\, =\, xz ~\  \text{ and }~\ \eta\,=\, y\, .
$$
Since $X$ commutes with $Y$, the action of $Y$ preserves the invariants of $X$.
We have
$$
Y(\xi)\,=\, 2\eta  ~\  \text{ and }~\ Y(\eta)\,=\, 1\, .
$$
So $Y$ on invariants of $X$ is
$$
Y\,=\, 2\eta\frac{\partial}{\partial \xi} +\frac{\partial}{\partial \eta}\, .
$$
Its invariants are given by $\frac{d\xi}{2\eta}\,=\, \frac{d\eta}{1}$. So the
basic invariant is
$$
\xi-\eta^2\,=\, xz-y^2\, .
$$

The next example is from \cite{Ib1}.

Take the following two vector fields on ${\mathbb R}^4$ with coordinates $(x\, ,y
\, ,z\, ,w)$:
$$
X_1\, :=\, (0\, ,z\, ,-y\, ,0)\, , X_2\, :=\, (1\, , w\, ,0\, , y)\, .
$$
Its row reduced echelon form $Y_1\, , Y_2$ is not closed under Lie bracket. We have
$$
[Y_1\, , Y_2]\, = \, X_3\, :=\, (0\, ,0\, ,-w/z\, ,-1)\, .
$$
The row reduced echelon form for $X_1\, ,X_2\, ,X_3$ is
$$
\begin{pmatrix}
1 & 0 & 0 & 0\\
0 & 1 & 0 & y/w\\
0 & 0 & 1 & z/w\\
\end{pmatrix}
$$
which give commutative vector fields. Consequently, the joint invariant
is $y^2+z^2-w^2$.

\section{More examples}

An efficient way to get invariants of a solvable algebra $L$ is to first determine
the joint invariants of the commutator algebra --- which is always nilpotent and thus one
can use the central series for systematic reductions --- and then find the joint invariants
of the full algebra as they are the same as those of $L/L'$ on the invariants of $L'$.

Also for semi-direct products $L \rtimes V$ one can first find the joint invariants of
$V$, and then the invariants of $L$ on the invariants of $V$ to find the joint
invariants of the full algebra.

Before giving examples, let us record the formulas for the
fundamental vector fields as differential operators in the adjoint and coadjoint
representations of Lie groups.

Let $\mathcal L$ be a finite dimensional Lie algebra, and let
$X_{1}\, ,\cdots\, ,X_{d}$ be a basis of $\mathcal L$.
Let $\omega_{1}\, ,\cdots\, ,\omega_{d}$ be the dual basis of ${\mathcal L}^*$.

For $X\,\in\, \mathcal L$, the fundamental vector fields $X_{\mathcal L}$ and $X_{{\mathcal
L}^*}$ corresponding to $X$ in the adjoint and coadjoint representations are given as
differential operators by the formulas:
$$
X_{\mathcal L}\,=\, \sum_{1\leq i, j \leq d} x_{i}\omega_{j}([X\, ,X_{i}])
\frac{\partial}{\partial{x_j}} ~\ \text{ and }~\
X_{{\mathcal L}^*} \,=\, -\sum_{1\leq i, j \leq d} x_{i}\omega_{i}([X\, ,X_{j}])
\frac{\partial}{\partial{x_j}}\, .
$$

Several examples of invariants of solvable algebras are computed in \cite{Nd, Pe1}.
Also invariants of real low dimensional algebras and some general classical algebras
are calculated in several papers, for example \cite{Pa1, Pa2, Pe2, PN}. We now give some
examples of fundamental invariants of certain solvable Lie algebras and Lie algebras
of low rank.

\subsection{Examples from \cite{BPP}}

For the convenience of the reader, we will refer to the online version of the
paper \cite{BPP} --- available at http://arxiv.org/pdf/math-ph/0602046.pdf.

\subsubsection{Example 1} We will use the variable $x\, ,y\, , z\, ,w$ for the variable
$\{e_i\}_{i=1}^4$ in Example 1 of \cite{BPP}.

After writing the matrix of the operators in the coadjoint representation, Maple
directly gives two joint invariant, one of which is in integral form. Working with the
reduced row echelon form we easily get one invariant
$$
I_1\, =\, (x^2+y^2)\exp (-2b\cdot \tan^{-1}(y/x))\, .
$$
A second invariant can be obtained by using elementary implications like
$$
\frac{a}{b}\,=\, \frac{c}{d} \, \Rightarrow\, \frac{a}{b}\,=\, \frac{\lambda a+\mu
c}{\lambda b +\mu d}\, .
$$
This gives a second independent invariant
$$
I_2\,=\, \frac{w^{2b}}{(x^2+y^2)^a}\, ;
$$
this corrects a misprint in this example from \cite{BPP}.

\subsubsection{Example 2} We will use the variable $s\, , w\, ,x\, ,y\, , z$ for the variable
$\{e_i\}_{i=1}^5$ in Example 2 of \cite{BPP}.

After writing down the matrix of coefficients of the operators in the coadjoint
representation corresponding to the given basis and using the operators corresponding
to the row reduced form, we find that there is only basic joint invariant
$$
\frac{w -s\cdot \ln s}{s}\, .
$$
Maple gives this directly --- without any row reductions.

\subsubsection{Example 3} We will use the variable $s\, , w\, ,x\, ,y\, , z$ for the variable
$\{e_i\}_{i=1}^5$ in Example 3 of \cite{BPP}.

Using the same procedure as in Example 2, Maple gives directly the invariant
$$\frac{xw+zs}{s}\, .$$

\subsubsection{Example 4} We will use the variable $r\, , s\, , w\, ,x\, ,y\, , z$ for
the variable $\{e_i\}_{i=1}^6$ in Example 4 of \cite{BPP}.

Maple cannot find directly joint invariants form the matrix of operators for the
coadjoint representation. However, when one works with the row reduced echelon form,
the situation simplifies dramatically. One gets two basic invariants
$$
I_1\,=\, r^{-2b}(x^2+w^2)\exp (-2a\cdot \tan^{-1}(w/x))\ ~ \text{ and }~ \
I_2\,=\, \frac{s}{r} - \frac{1}{2a}\ln \frac{x^2+w^2}{r^{2b}}\, .
$$

\subsection{Invariants of $\text{sl}(3,{\mathbb R})$ in its adjoint and coadjoint
representations}

The non-zero commutation relations are
\begin{small}
$$
[e_1, e_2] = e_2, [e_1, e_3] = 2 e_3, [e_1, e_4] = -e_4, [e_1, e_6] = e_6, [e_1, e_7] = -2 e_7, [e_1, e_8] = -e_8,
[e_2, e_4] = e_1-e_5,
$$
$$
[e_2, e_5] = e_2, [e_2, e_6] = e_3, [e_2, e_7] = -e_8, [e_3, e_4] = -e_6, [e_3, e_5] = -e_3, [e_3, e_7] = e_1, [e_3, e_8] = e_2, [e_4, e_5] = -e_4,
$$
$$
 [e_4, e_8] = -e_7, [e_5, e_6] = 2 e_6, [e_5, e_7] = -e_7, [e_5, e_8] = -2 e_8, [e_6, e_7] = e_4, [e_6, e_8] = e_5
$$
\end{small}
Writing the operators $\displaystyle \sum_{i=1}^{8} x_i X_i$ as $[x_1,x_2,\cdots,x_8]$, the coadjoint representation of the basis of $\text{sl}(3,{\mathbb R})$ is \begin{eqnarray*}
X_1 &=& [0,-x_{{2}},-2\,x_{{3}},x_{{4}},0,-x_{{6}},2\,x_{{7}},x_{{8}}]\\
X_2 &=&[x_{{2}},0,0,x_{{5}}-x_{{1}},-x_{{2}},-x_{{3}},x_{{8}},0] \\
X_3 &=& [2\,x_{{3}},0,0,x_{{6}},x_{{3}},0,-x_{{1}},-x_{{2}}]   \\
X_4 &=& [-x_{{4}},-x_{{5}}+x_{{1}},-x_{{6}},0,x_{{4}},0,0,x_{{7}}]  \\
X_5 &=& [0,x_{{2}},-x_{{3}},-x_{{4}},0,-2\,x_{{6}},x_{{7}},2\,x_{{8}}]  \\
X_6 &=&  [x_{{6}},x_{{3}},0,0,2\,x_{{6}},0,-x_{{4}},-x_{{5}}] \\
X_7 &=& [-2\,x_{{7}},-x_{{8}},x_{{1}},0,-x_{{7}},x_{{4}},0,0]  \\
X_8 &=&  [-x_{{8}},0,x_{{2}},-x_{{7}},-2\,x_{{8}},x_{{5}},0,0] \\
\end{eqnarray*}
The reduced echelon form
\begin{equation}
\begin{tiny}
 \left[ \begin {array}{cccccccc}
 1&0&0&0&0&0&{\frac {2 {x_{{6}
}}^{2}x_{{8}}-x_{{3}}x_{{7}}x_{{6}}-x_{{2}}x_{{4}}x_{{6}}-x_{{3}}x_{{4
}}x_{{5}}+2 x_{{1}}x_{{4}}x_{{3}}+2 x_{{1}}x_{{5}}x_{{6}}-2 {x_{{1}
}}^{2}x_{{6}}}{3(-x_{{3}}x_{{5}}x_{{6}}+x_{{6}}x_{{3}}x_{{1}}-x_{{4}}{x_
{{3}}}^{2}+{x_{{6}}}^{2}x_{{2}})}}&{\frac {2 x_{{3}}x_{{8}}x_{{6
}}-{x_{{3}}}^{2}x_{{7}}-x_{{5}}x_{{6}}x_{{2}}+x_{{3}}{x_{{5}}}^{2}-x_{
{1}}x_{{3}}x_{{5}}-x_{{2}}x_{{4}}x_{{3}}+2 x_{{1}}x_{{2}}x_{{6}}}{3(x_
{{3}}x_{{5}}x_{{6}}-x_{{6}}x_{{3}}x_{{1}}+x_{{4}}{x_{{3}}}^{2}-{x_{{6}
}}^{2}x_{{2}})}}\\
\noalign{\medskip}0&1&0&0&0&0&-{\frac {x_{{7}}{x_{{6
}}}^{2}+x_{{4}}x_{{6}}x_{{1}}-{x_{{4}}}^{2}x_{{3}}}{-x_{{3}}x_{{5}}x_{
{6}}+x_{{6}}x_{{3}}x_{{1}}-x_{{4}}{x_{{3}}}^{2}+{x_{{6}}}^{2}x_{{2}}}}
&-{\frac {-x_{{3}}x_{{7}}x_{{6}}-x_{{3}}x_{{4}}x_{{5}}+x_{{2}}x_{{4}}x
_{{6}}}{-x_{{3}}x_{{5}}x_{{6}}+x_{{6}}x_{{3}}x_{{1}}-x_{{4}}{x_{{3}}}^
{2}+{x_{{6}}}^{2}x_{{2}}}}\\
\noalign{\medskip}0&0&1&0&0&0&-{\frac {x_
{{7}}x_{{6}}x_{{1}}-x_{{7}}x_{{4}}x_{{3}}-x_{{7}}x_{{5}}x_{{6}}+x_{{6}
}x_{{4}}x_{{8}}}{-x_{{3}}x_{{5}}x_{{6}}+x_{{6}}x_{{3}}x_{{1}}-x_{{4}}{
x_{{3}}}^{2}+{x_{{6}}}^{2}x_{{2}}}}&-{\frac {x_{{7}}x_{{6}}x_{{2}}-x_{
{4}}x_{{3}}x_{{8}}}{-x_{{3}}x_{{5}}x_{{6}}+x_{{6}}x_{{3}}x_{{1}}-x_{{4
}}{x_{{3}}}^{2}+{x_{{6}}}^{2}x_{{2}}}}\\
\noalign{\medskip}0&0&0&1&0&0
&-{\frac {x_{{1}}x_{{2}}x_{{6}}+x_{{3}}x_{{8}}x_{{6}}-x_{{2}}x_{{4}}x_
{{3}}}{-x_{{3}}x_{{5}}x_{{6}}+x_{{6}}x_{{3}}x_{{1}}-x_{{4}}{x_{{3}}}^{
2}+{x_{{6}}}^{2}x_{{2}}}}&-{\frac {-{x_{{3}}}^{2}x_{{8}}-x_{{2}}x_{{5}
}x_{{3}}+{x_{{2}}}^{2}x_{{6}}}{-x_{{3}}x_{{5}}x_{{6}}+x_{{6}}x_{{3}}x_
{{1}}-x_{{4}}{x_{{3}}}^{2}+{x_{{6}}}^{2}x_{{2}}}}\\
\noalign{\medskip}0
&0&0&0&1&0&{\frac {-x_{{2}}x_{{4}}x_{{6}}+2 x_{{3}}x_{{4}}x_{{5}
}-x_{{1}}x_{{4}}x_{{3}}-{x_{{6}}}^{2}x_{{8}}+2 x_{{3}}x_{{7}}x_{{6}}-
x_{{1}}x_{{5}}x_{{6}}+{x_{{1}}}^{2}x_{{6}}}{3(-x_{{3}}x_{{5}}x_{{6}}+x_{
{6}}x_{{3}}x_{{1}}-x_{{4}}{x_{{3}}}^{2}+{x_{{6}}}^{2}x_{{2}})}}&{
\frac {-2 x_{{5}}x_{{6}}x_{{2}}+2 x_{{3}}{x_{{5}}}^{2}-2 x_{{1}}x_{
{3}}x_{{5}}+x_{{3}}x_{{8}}x_{{6}}-2 {x_{{3}}}^{2}x_{{7}}+x_{{2}}x_{{4
}}x_{{3}}+x_{{1}}x_{{2}}x_{{6}}}{3(-x_{{3}}x_{{5}}x_{{6}}+x_{{6}}x_{{3}}
x_{{1}}-x_{{4}}{x_{{3}}}^{2}+{x_{{6}}}^{2}x_{{2}})}}
\\
\noalign{\medskip}0&0&0&0&0&1&-{\frac {x_{{7}}x_{{6}}x_{{2}}-x_{{4}
}x_{{3}}x_{{8}}}{-x_{{3}}x_{{5}}x_{{6}}+x_{{6}}x_{{3}}x_{{1}}-x_{{4}}{
x_{{3}}}^{2}+{x_{{6}}}^{2}x_{{2}}}}&-{\frac {x_{{8}}x_{{6}}x_{{2}}-x_{
{8}}x_{{5}}x_{{3}}+x_{{8}}x_{{3}}x_{{1}}-x_{{3}}x_{{7}}x_{{2}}}{-x_{{3
}}x_{{5}}x_{{6}}+x_{{6}}x_{{3}}x_{{1}}-x_{{4}}{x_{{3}}}^{2}+{x_{{6}}}^
{2}x_{{2}}}}\\
\noalign{\medskip}0&0&0&0&0&0&0&0\\
\noalign{\medskip}0
&0&0&0&0&0&0&0 \nonumber
\end {array} \right]
\end{tiny}
\end{equation}
leads to commuting operators, and implies that there are two joint invariants which can be found using Maple as
\begin{tiny}
\begin{eqnarray*}
I_1 &=&  {x_{{5}}}^{2}+{x_{{1}}}^{2}-x_{{1}}x_{{5}}+3 x_{{7}}x_{{3}}+3 x_{{8}
}x_{{6}}+3 x_{{2}}x_{{4}} \\
I_2 &=& 2 {x_{{1}}}^{3}-3 x_{{5}}{x_{{1}}}^{2}+9 x_{{2}}x_{{4}}x_{{1}}-3 x
_{{1}}{x_{{5}}}^{2}-18 x_{{1}}x_{{8}}x_{{6}}+9 x_{{7}}x_{{3}}x_{{1}}
+2 {x_{{5}}}^{3}+9 x_{{5}}x_{{8}}x_{{6}}-18 x_{{7}}x_{{5}}x_{{3}}+9
 x_{{5}}x_{{2}}x_{{4}}+27 x_{{7}}x_{{6}}x_{{2}}+27 x_{{4}}x_{{3}}x_
{{8}}
\end{eqnarray*}
\end{tiny}

The adjoint representation of the basis of $\text{sl}(3,{\mathbb R})$ is
\begin{eqnarray*}
X_1 &=& [0,x_{{2}},2\,x_{{3}},-x_{{4}},0,x_{{6}},-2\,x_{{7}},-x_{{8}}]\\
X_2 &=& [x_{{4}},x_{{5}}-x_{{1}},x_{{6}},0,-x_{{4}},0,0,-x_{{7}}] \\
X_3 &=& [x_{{7}},x_{{8}},-x_{{5}}-2\,x_{{1}},0,0,-x_{{4}},0,0]   \\
X_4 &=&  [-x_{{2}},0,0,-x_{{5}}+x_{{1}},x_{{2}},x_{{3}},-x_{{8}},0] \\
X_5 &=&  [0,-x_{{2}},x_{{3}},x_{{4}},0,2\,x_{{6}},-x_{{7}},-2\,x_{{8}}] \\
X_6 &=&  [0,0,-x_{{2}},x_{{7}},x_{{8}},-2\,x_{{5}}-x_{{1}},0,0] \\
X_7 &=&  [-x_{{3}},0,0,-x_{{6}},0,0,x_{{5}}+2\,x_{{1}},x_{{2}}] \\
X_8 &=& [0,-x_{{3}},0,0,-x_{{6}},0,x_{{4}},2\,x_{{5}}+x_{{1}}] \\
\end{eqnarray*}
The reduced echelon form
\begin{equation}
\begin{tiny}
\left[ \begin {array}{ccccccccc}
1&0&0&0&0&0&-{\frac {-2\,x_{{3}}x_{{
1}}x_{{4}}-x_{{3}}x_{{5}}x_{{4}}+x_{{4}}x_{{2}}x_{{6}}-{x_{{6}}}^{2}x_
{{8}}-x_{{6}}{x_{{5}}}^{2}-x_{{6}}x_{{5}}x_{{1}}+2\,x_{{6}}{x_{{1}}}^{
2}}{x_{{3}}x_{{1}}x_{{6}}+{x_{{6}}}^{2}x_{{2}}-{x_{{3}}}^{2}x_{{4}}-x_
{{5}}x_{{3}}x_{{6}}}}&-{\frac {2\,x_{{2}}x_{{6}}x_{{1}}-x_{{3}}x_{{4}}
x_{{2}}+x_{{6}}x_{{8}}x_{{3}}+x_{{6}}x_{{5}}x_{{2}}}{x_{{3}}x_{{1}}x_{
{6}}+{x_{{6}}}^{2}x_{{2}}-{x_{{3}}}^{2}x_{{4}}-x_{{5}}x_{{3}}x_{{6}}}}\\
\noalign{\medskip}0&1&0&0&0&0&-{\frac {x_{{4}}x_{{5}}x_{{6}}+2\,x
_{{4}}x_{{6}}x_{{1}}+{x_{{6}}}^{2}x_{{7}}-{x_{{4}}}^{2}x_{{3}}}{x_{{3}
}x_{{1}}x_{{6}}+{x_{{6}}}^{2}x_{{2}}-{x_{{3}}}^{2}x_{{4}}-x_{{5}}x_{{3
}}x_{{6}}}}&{\frac {x_{{3}}x_{{1}}x_{{4}}+x_{{7}}x_{{3}}x_{{6}}+2\,x_{
{3}}x_{{5}}x_{{4}}-x_{{4}}x_{{2}}x_{{6}}}{x_{{3}}x_{{1}}x_{{6}}+{x_{{6
}}}^{2}x_{{2}}-{x_{{3}}}^{2}x_{{4}}-x_{{5}}x_{{3}}x_{{6}}}}
\\
\noalign{\medskip}0&0&1&0&0&0&-{\frac {x_{{6}}x_{{7}}x_{{1}}+x_{{6}
}x_{{4}}x_{{8}}-x_{{7}}x_{{4}}x_{{3}}-x_{{6}}x_{{7}}x_{{5}}}{x_{{3}}x_
{{1}}x_{{6}}+{x_{{6}}}^{2}x_{{2}}-{x_{{3}}}^{2}x_{{4}}-x_{{5}}x_{{3}}x
_{{6}}}}&-{\frac {x_{{2}}x_{{6}}x_{{7}}-x_{{4}}x_{{8}}x_{{3}}}{x_{{3}}
x_{{1}}x_{{6}}+{x_{{6}}}^{2}x_{{2}}-{x_{{3}}}^{2}x_{{4}}-x_{{5}}x_{{3}
}x_{{6}}}}\\
\noalign{\medskip}0&0&0&1&0&0&-{\frac {2\,x_{{2}}x_{{6}
}x_{{1}}-x_{{3}}x_{{4}}x_{{2}}+x_{{6}}x_{{8}}x_{{3}}+x_{{6}}x_{{5}}x_{
{2}}}{x_{{3}}x_{{1}}x_{{6}}+{x_{{6}}}^{2}x_{{2}}-{x_{{3}}}^{2}x_{{4}}-
x_{{5}}x_{{3}}x_{{6}}}}&{\frac {x_{{2}}x_{{3}}x_{{1}}+{x_{{3}}}^{2}x_{
{8}}+2\,x_{{3}}x_{{5}}x_{{2}}-{x_{{2}}}^{2}x_{{6}}}{x_{{3}}x_{{1}}x_{{
6}}+{x_{{6}}}^{2}x_{{2}}-{x_{{3}}}^{2}x_{{4}}-x_{{5}}x_{{3}}x_{{6}}}}
\\
\noalign{\medskip}0&0&0&0&1&0&{\frac {x_{{3}}x_{{1}}x_{{4}}+x_{{7}}
x_{{3}}x_{{6}}+2\,x_{{3}}x_{{5}}x_{{4}}-x_{{4}}x_{{2}}x_{{6}}}{x_{{3}}
x_{{1}}x_{{6}}+{x_{{6}}}^{2}x_{{2}}-{x_{{3}}}^{2}x_{{4}}-x_{{5}}x_{{3}
}x_{{6}}}}&-{\frac {x_{{3}}{x_{{1}}}^{2}+x_{{2}}x_{{6}}x_{{1}}+x_{{3}}
x_{{5}}x_{{1}}+{x_{{3}}}^{2}x_{{7}}-2\,x_{{3}}{x_{{5}}}^{2}-x_{{3}}x_{
{4}}x_{{2}}+2\,x_{{6}}x_{{5}}x_{{2}}}{x_{{3}}x_{{1}}x_{{6}}+{x_{{6}}}^
{2}x_{{2}}-{x_{{3}}}^{2}x_{{4}}-x_{{5}}x_{{3}}x_{{6}}}}
\\
\noalign{\medskip}0&0&0&0&0&1&-{\frac {x_{{2}}x_{{6}}x_{{7}}-x_{{4}
}x_{{8}}x_{{3}}}{x_{{3}}x_{{1}}x_{{6}}+{x_{{6}}}^{2}x_{{2}}-{x_{{3}}}^
{2}x_{{4}}-x_{{5}}x_{{3}}x_{{6}}}}&-{\frac {x_{{3}}x_{{8}}x_{{1}}+x_{{
2}}x_{{8}}x_{{6}}-x_{{2}}x_{{7}}x_{{3}}-x_{{3}}x_{{8}}x_{{5}}}{x_{{3}}
x_{{1}}x_{{6}}+{x_{{6}}}^{2}x_{{2}}-{x_{{3}}}^{2}x_{{4}}-x_{{5}}x_{{3}
}x_{{6}}}}\\
\noalign{\medskip}0&0&0&0&0&0&0&0
\\
\noalign{\medskip}0&0&0&0&0&0&0&0\end {array} \right] \nonumber
\end{tiny}
\end{equation}
leads to commuting operators, and implies that there are two joint invariants which can be found using Maple as
\begin{eqnarray*}
I_1 &=& {x_{{5}}}^{2}+x_{{1}}x_{{5}}+{x_{{1}}}^{2}+x_{{7}}x_{{3}}+x_{{8}}x_{{6
}}+x_{{4}}x_{{2}}
  \\
I_2 &=& -{x_{{1}}}^{2}x_{{5}}-x_{{1}}x_{{6}}x_{{8}}+x_{{1}}x_{{4}}x_{{2}}-x_{{
1}}{x_{{5}}}^{2}-x_{{3}}x_{{7}}x_{{5}}+x_{{4}}x_{{8}}x_{{3}}+x_{{2}}x_
{{6}}x_{{7}}+x_{{4}}x_{{5}}x_{{2}}
\end{eqnarray*}

\subsection{Invariants of forms of $\text{so}(4)$ in their adjoint and coadjoint
representations}  The basic invariants for real forms of
$\text{so}(4)$ in suitable coordinates obtained as in 3.2 are \\

$\text{so}(4): $\\ [0.1in]
$x^2_4+x^2_3+2x_4x_3+(x_1+x_6)^2+(x_5-x_2)^2$\\
$x^2_4+x^2_3-2x_4x_3+(x_1-x_6)^2+(x_5+x_2)^2$\\

$\text{so}(2,2): $\\ [0.1in]
$x^2_4+x^2_3-2x_4x_3+(x_5+x_2)^2-(x_1-x_6)^2$\\
$x^2_4+x^2_3+2x_4x_3+(x_2-x_5)^2-(x_1+x_6)^2$\\

$\text{so}(1,3):$\\[0.1in]
$-x^2_4+x^2_3-2Ix_4x_3-(x_1+Ix_6)^2+(x_5+Ix_2)^2$\\
$-x^2_4+x^2_3+2Ix_4x_3-(x_1-Ix_6)^2+(x_5-Ix_2)^2$\\

The real invariants are obtained from taking the real and imaginary
parts of either of the above two invariants.

\subsection{Concluding remarks}

The commuting vector fields which give the invariants of the exceptional groups can
also be computed because explicit structure constants, which are
programmable, are available -- as indicated e.g in \cite{Az1}, \cite[p. 9]{Az2};
see also \cite{Va}.

The exceptional groups are also of interest to theoretical physicists \cite{Ca}, \cite{Ra}.
In certain cases, the joint invariants in the fundamental representations of certain
exceptional groups can also be obtained algorithmically.

For example, one can realize $D_{4}$ is the Levi complement of a maximal parabolic
subgroups of $D_{5}$ as in \cite{ABS}, use a choice of structure constants which are integers
and use triality to obtain $G_{2}$ as a subgroup of $D_{5}$ with its maximal torus as a
subgroup of a maximal torus of $D_{5}$. Then the root vector corresponding to the simple
root of $D_{5}$ which is not a simple root of $D_{4}$ would be a high weight vector for
$G_{2}$ and it translates under $G_{2}$ would give the seven dimensional fundamental
representation of $G_{2}$.

\section*{Acknowledgements}

We thank K.-H. Neeb for a very helpful correspondence.

\end{document}